\documentclass[%
a4paper,							
11pt,								
bibliography=totoc,						
abstracton,					
]
{scrartcl}
\usepackage[a4paper,left=3.0cm,right=2.5cm, top=2.5cm, bottom=3.0cm]{geometry}
\usepackage[headsepline=.4pt,footsepline=.4pt,automark,autooneside=false,]{scrlayer-scrpage}
\clearpairofpagestyles
\automark[subsection]{section} 
\automark*[subsubsection]{subsection}
\ihead{\scriptsize\rightmark} 
\usepackage[dvipsnames]{xcolor} 
\usepackage[colorinlistoftodos,prependcaption]{todonotes}

\usepackage[ngerman, english]{babel}
\usepackage[T1]{fontenc}
\usepackage[utf8]{inputenc}
\usepackage{lmodern} 
\usepackage{csquotes}
\usepackage[section]{placeins} 
\usepackage{calc}
\usepackage{xspace}

\usepackage{amsmath}
\usepackage{amssymb}
\usepackage{amsfonts}
\usepackage{amsthm,thmtools}
\usepackage{mathtools}
\usepackage{bbm}

\theoremstyle{plain}

\theoremstyle{definition}

\theoremstyle{remark}

\usepackage{graphicx,float}
\usepackage[format=plain]{caption}

\usepackage{setspace}

\usepackage{xcolor}

\usepackage{fancyvrb}
\usepackage{xcolor}
\usepackage{paralist}

\usepackage{tabularx}
\usepackage{ragged2e}
\newcolumntype{C}{>{\Centering}X}
\usepackage{booktabs}
\usepackage{multirow}
\usepackage{multicol}
\usepackage{pdflscape}
\usepackage{diagbox}

\usepackage{tikz}
\usetikzlibrary{calc,positioning,arrows.meta}

\usepackage{todonotes}
\usepackage{makecell}
\usepackage[
	backend=bibtex,
	style=authoryear-comp,
	maxbibnames=9,
	maxcitenames=2,
	dashed=false,
	natbib=true,
	sortcites=true,
	block=space
]{biblatex}

\usepackage{orcidlink}
\newcommand{\orcid}[1]{\textsuperscript{\orcidlink{#1}}}
\usepackage{titleref}
\usepackage{hyperref}
\usepackage{cleveref}
\hypersetup{
    colorlinks,
    linkcolor={red},
    citecolor={blue},
    urlcolor={blue}
}

\newenvironment{key}[1][]{ {\noindent \bf Keywords#1: } \rmfamily
}{}
\newenvironment{code}[1][]{ {\noindent \bf Code availability#1: } \rmfamily
}{}%
\newenvironment{acknowledgements}[1][]{ {\noindent \bf Acknowledgements#1: } \rmfamily
}{}%

\newcommand{\Q}{\mathbb{Q}}
\newcommand{\E}{\mathbb{E}}
\newcommand{\R}{\mathbb{R}}

\newcommand{\abs}[1]{\lvert #1 \rvert}

\DeclareMathOperator*{\argmax}{argmax}

\SaveVerb{python}=(Intel-)Python 3.10=

\SaveVerb{tensorflow}=Tensorflow 2.8.0=

\SaveVerb{softmax}=softmax=

\SaveVerb{relu}=relu=

\SaveVerb{matlab}=Matlab 2023a=

\SaveVerb{ga}=ga=

\SaveVerb{fmincon}=fmincon=

\SaveVerb{fitdist}=fitdist=

\SaveVerb{lsqnonlin}=lsqnonlin=

\definecolor{comments1}{rgb}{0.858, 0.0, 0.0}

\begin{document}

\begin{titlepage}
\begin{center}
    {\huge \bfseries \sffamily Joint Stochastic Optimal Control and Stopping in Aquaculture: Finite-Difference and PINN-Based Approaches}
\end{center}
\renewcommand{\thefootnote}{\fnsymbol{footnote}}
\footnotetext[1]{\href{mailto:kevin.kamm@umu.se}{kevin.kamm@umu.se}, Department of Mathematics and Statistics, Umeå University, Sweden.}
\vspace*{-.33cm}
\begin{center}
    \Large
    Kevin Kamm%
    \footnotemark[1]{}\textsuperscript{,}%
    \orcid{0000-0003-2881-0905}
\end{center}
\begin{center}
    \large
    \today
\end{center}
\vspace*{.5cm}

\renewcommand{\thefootnote}{\arabic{footnote}}

\begin{abstract}
This paper studies a joint stochastic optimal control and stopping (JCtrlOS) problem motivated by aquaculture operations, where the objective is to  maximize farm profit through an optimal feeding strategy and harvesting time under stochastic price dynamics. We introduce a simplified aquaculture model capturing essential biological and economic features, distinguishing between biologically optimal and economically optimal feeding strategies. The problem is formulated as a Hamilton-Jacobi-Bellman variational inequality and corresponding free boundary problem. We develop two numerical solution approaches: First, a finite difference scheme that serves as a benchmark, and second a Physics-Informed Neural Network (PINN)-based method, combined with a deep optimal stopping (DeepOS) algorithm to improve stopping time accuracy. Numerical experiments demonstrate that while finite differences perform well in medium-dimensional settings, the PINN approach achieves comparable accuracy and is more scalable to higher dimensions where grid-based methods become infeasible. The results confirm that jointly optimizing feeding and harvesting decisions outperforms strategies that neglect either control or stopping.
\end{abstract}

\medskip
\begin{key}
Stochastic Optimal Control, Optimal Stopping, Deep Learning, PINN, Finite Differences, HJB, Numerics, Aquaculture
\end{key}

\medskip
\begin{code}
The Python code and data sets to produce the numerical experiments are available at 
\url{https://github.com/kevinkamm/JointOptimalCtrlStopping_Aquaculture}.\\[-\baselineskip]
\end{code}

\medskip
\begin{acknowledgements}
This research has been part of the AquaRisk project, project number 352447, funded by the Norwegian Research Council.
\end{acknowledgements}
\vfill
\begin{center}
    {\small \textbf{Graphical Abstract}}\\\noindent
    \includegraphics[width=.75\columnwidth]{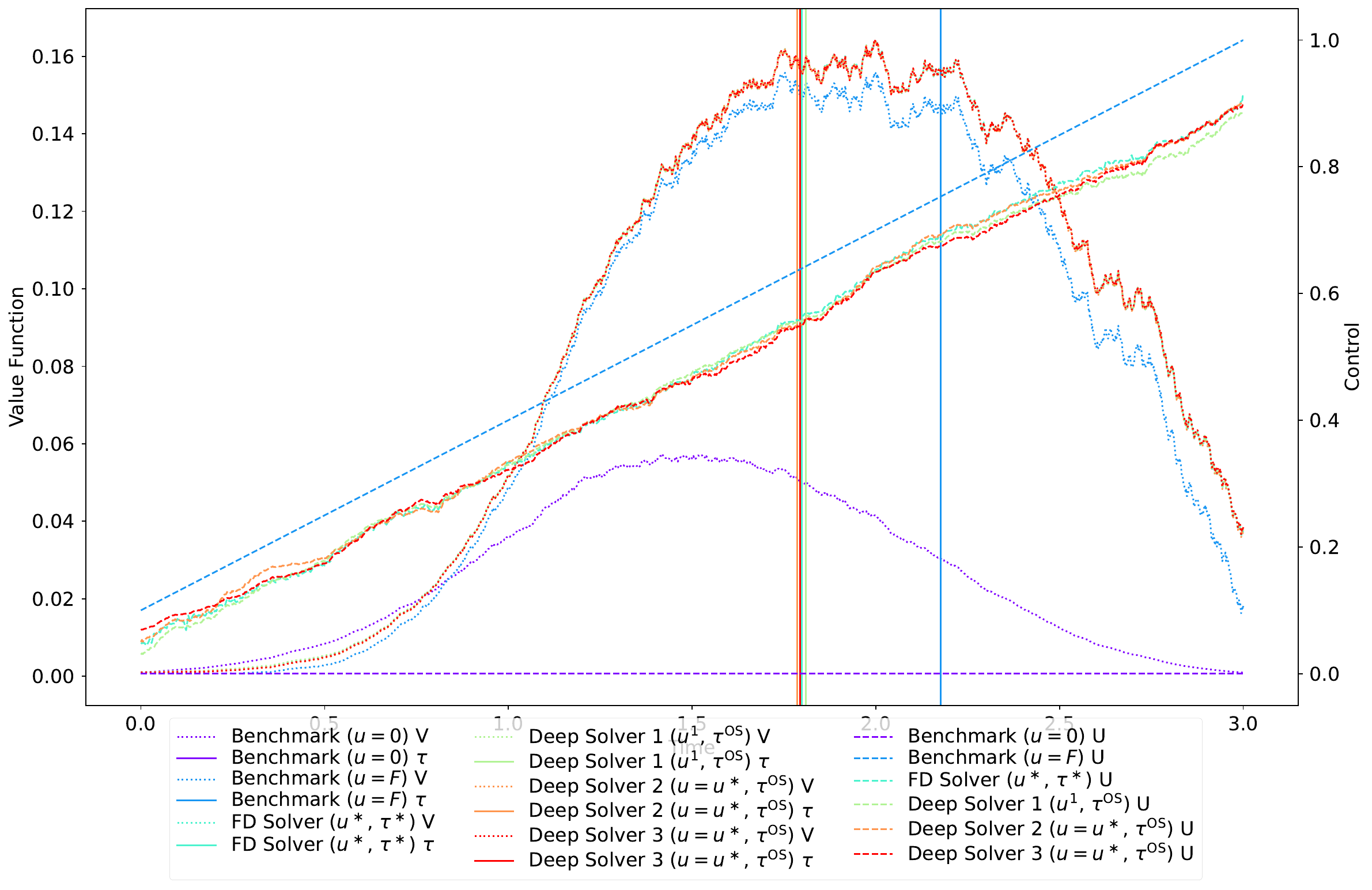}
\end{center}

\thispagestyle{empty}
\end{titlepage}

\newpage
\pagestyle{scrheadings}
\ihead{\scriptsize\rightmark}
\pagenumbering{arabic}

\section{Introduction}\label{sec:Introduction}
In this paper, we will investigate a joint stochastic optimal control and stopping problem (JCtrlOS) with medium dimensionality\footnote{With medium dimensionality we mean that classical methods can still be applied, which is the case for a five-dimensional setting as in this paper.}. A JCtrlOS (see \eqref{eq:JCtrlOS}) is a mixture or combination of an optimal control problem and an optimal stopping problem. The overall goal is to find an optimal control strategy and an optimal stopping time such that a certain objective functional is maximized. Such problems arise naturally in finance, e.g., portfolio optimization with American options or resource management. We will consider an application in aquaculture, where the goal is to find an optimal feeding strategy and an optimal harvesting time such that the profit of the farm is maximized. However, it is not the goal of this paper to provide a realistic model for aquaculture, but rather to provide a framework that can be used in practice and to compare two different numerical methods to solve the problem. Therefore, we pick a toy model (see \Cref{sec:Setting}) for the aquaculture dynamics that is simple enough to be solved by a finite difference scheme and still captures the main features of the problem. Then, we will compare it to a Deep Learning approach based on Physics-Informed Neural Networks (PINNs) to solve the same problem and combine it with the DeepOS method of \citet{Becker2020} to improve the solution. We will see that the PINN approach is able to solve the problem with similar accuracy as the finite difference scheme and is therefore a viable alternative. We conjecture that for higher dimensional problems, where finite difference schemes are not applicable anymore, the PINN approach will be a good alternative to solve JCtrlOS problems.

A notable feature of the toy model is that the stochastic optimal control can only act on the biological part of the model, i.e., the weight and population dynamics, which are assumed to be deterministic, but not on the price dynamics. This is motivated by the fact that from a biological models from theoretical studies or mean behaviour from empirical studies can be used to determine biological parameters, such as intrinsic mortality and growth of the species or a feeding strategy that maximizes growth. This feeding strategy will be called in this paper \emph{biological optimal feeding strategy} and denoted by $f_t$. However, the farm operator can deviate from this biological optimal feeding strategy in order to maximize the profit of the farm by considering stochastic price dynamics. This feeding strategy will be called in this paper \emph{economic optimal feeding strategy}. The goal is to study whether or not it is beneficial to deviate from the biological optimal feeding strategy in order to maximize the profit of the farm by considering stochastic price dynamics. And the answer will be affirmative in the numerical experiments in our toy model. 
\paragraph*{Literature review.}\label{sec:LiteratureReview}
Let us give a brief overview of the relevant literature to this paper and recent developments in the respective fields.

\subparagraph*{Aquaculture Modelling.}
For an overview of aquaculture modelling and recent developments, we refer the reader to \citet{EK2024,EK2025} and references therein. We reduced the complexity of the models in \citet{EK2024,EK2025} significantly to be able to solve the problem with a finite difference scheme and still capture the main features of the problem.

\subparagraph*{Joint Optimal Control and Stopping Problems.} In the literature, joint optimal control and stopping is also referred to as mixed optimal control and stopping by \citet{Reisinger2019} and combined optimal control and stopping by \citet{Dumitrescu2016}. Both papers treat the problem in a great generality. \citet{Dumitrescu2016} lay the theoretical foundation for existence and uniqueness in the case of non-linear expectations and \citet{Reisinger2019} provide a numerical scheme based on finite differences and prove convergence of the scheme. We will discuss alternative numerical schemes in this paper, which are based on solving the HJB variational inequality (HJB-VI) or the corresponding free boundary problem (HJB-FB) using PINNs.

\subparagraph*{Classical Tools for Optimal Control and Stopping Problems.}
There is a vast literature on classical methods for solving optimal control and stopping problems, including dynamic programming, the Hamilton-Jacobi-Bellman (HJB) equation, and variational methods. However, the combination of both control and stopping is less explored. The most influential work for this paper is the work by \citet{Kushner2001}. We adopt the methods presented in this book to solve the HJB-VI using a finite difference scheme. 

\subparagraph*{Deep Learning Tools for Optimal Control and Stopping Problems.}
There are many different streams of research that use Deep Learning to solve optimal control and stopping problems. The most relevant ones for this paper are the Deep Galerkin Method (DGM) by \citet{Sirignano2018} and the Physics-Informed Neural Networks (PINNs) by \citet{Raissi2019}. Both methods use neural networks to approximate the solution of partial differential equations (PDEs). 

We will also use ideas from \citet{Reppen2025}, who study deep learning methods for determining the free boundary in optimal stopping problems. We will not use exactly the same method, but noticed that the idea of a fuzzy boundary improved the quality of training the neural network significantly.

\medskip
To the best of our knowledge, this is the first paper to use Deep Learning methods to tackle high-dimensional joint optimal control and stopping problems. Especially, in aquaculture modelling, this is a novel approach, where either optimal control or optimal stopping is considered separately or the models are deterministic. The methodology and the deep learning approach can be easily generalized to other fields and applications, where joint optimal control and stopping problems arise.

\bigskip
In \Cref{sec:Setting}, we will introduce the toy model for aquaculture and formulate the joint optimal control and stopping problem. In \Cref{sec:HJBVI}, we will derive the HJB variational inequality and the corresponding free boundary problem. In \Cref{sec:NumExp}, we will present numerical experiments. Thus section is split in two major parts. In \Cref{sec:FDS}, we will show a finite difference scheme and in \Cref{sec:DS} a PINN approach to solve the problem. Finally, in \Cref{sec:Conclusion}, we will conclude and give an outlook on future research.

\section{The Model Framework}\label{sec:Setting}
We will distinguish between two types of optimality regarding feeding strategies in aquaculture. The first one is the optimal feeding strategy $f$ in the sense of maximal growth, i.e., the feeding strategy that maximizes the biomass of the fish. We will call this the \emph{biological optimal feeding rate}. The second one is the optimal feeding strategy $u$ in the sense of profit, i.e., the feeding strategy that maximizes the profit of the farm. We will call this the \emph{economic optimal feeding rate}.

We will assume that there exists an optimal (in the sense of maximal growth) deterministic feeding strategy $f\colon\R_+\to\R_+$ only dependent on time. This can be obtained from biological experiments or expert knowledge by considering metabolic rates for the particular fish species. We expect that a farm operator will have have suitable experience and knowledge to determine such a feeding strategy.

The goal is to study whether or not it is beneficial to deviate from this biological optimal feeding strategy $f$ in order to maximize the profit of the farm by considering stochastic price dynamics. We will denote this economic optimal feeding strategy by $u$.

To this end, we will consider the following toy farm dynamics for fixed biological optimal feeding strategy $f$ and determine the economic optimal feeding strategy $u$.

\paragraph*{Weight dynamics.}\label{sec:WeightDynamics}
We will denote for given biological optimal feeding rate $f$ and stochastic control $u$ the weight of the fish at time $t\geq 0$ by $w_t^{f,u}$ satisfying the following logistical growth dynamics
\begin{align*}
    d w_t^{f,u} &= \left(\gamma - \gamma^F \left(f_t-u_t\right)^2\right) \left( 1 - \left(\frac{w_t^{f,u}}{w_\infty}\right)^\nu \right) dt \eqqcolon b^w(t,w_t^{f,u};f_t,u_t) dt, \\
    w_0^{f,u} &= w_0>0, \\
\end{align*}
where $\gamma>0$ is the intrinsic growth rate, which is optimal when feeding the biological optimal feeding rate $f_t$, i.e. $u_t=f_t$, and $\gamma^F>0$ is the growth deficit due to under- or overfeeding.\footnote{We choose a quadratic penalty term for simplicity to penalize over- and underfeeding equally. In an enclosed farm overfeeding can reduce the water quality and therefore be harmful to the population and underfeeding can lead to malnutrition. In the optimization problem later on, overfeeding will be penalized by the additional feeding costs and therefore never be optimal in our setting.}
The parameter $w_\infty>0$ is the asymptotic weight of the fish, i.e., the weight that the fish will reach if fed optimally for an infinite amount of time. The initial weight of the fish is given by $w_0>0$ and the parameter $\nu$ determines the speed of convergence to the asymptotic weight $w_\infty$. For $\nu=1$ we obtain the classical logistical growth model. For $\nu>1$ the fish will reach the asymptotic weight faster in the beginning and slower at the end, while for $\nu<1$ it is the other way around.


\paragraph*{Population dynamics.}\label{sec:PopulationDynamics}
We will denote by $h_t^{f,u}$ the number of fish at time $t\geq 0$ when applying the feeding strategy $f_t$ and the economic optimal feeding strategy $u_t$. We assume that the dynamics of $h_t^{f,u}$ is given by the following exponential growth model assuming that there is no reproduction of the fish in the farm
\begin{align*}
    d h_t^{f,u} &= \left(-\mu - \mu^F \left(f_t-u_t\right)^2\right) h_t^{f,u} dt \eqqcolon b^h(t,h_t^{f,u};f_t,u_t)dt , \\
    h_0^{f,u} &= h_0>0, \\
\end{align*}
where $\mu>0$ is the intrinsic mortality rate of the fish and $\mu^F>0$ the mortality increase due to under- or overfeeding. The initial number of fish in the farm is given by $h_0>0$.


\paragraph*{Price dynamics.}\label{sec:PriceDynamics}
We will denote by $P^F_t$ and $P^B_t$ the price of feed and the price of biomass at time $t\geq 0$. We will assume that the price dynamics are both by two geometric Brownian motions with dynamics under a risk-neutral measure $\Q$ given by
\begin{align*}
    d P^F_t &= r P^F_t dt + \sigma_F P^F_t dW^F_t \eqqcolon b^{p^F}(t,P^F_t)dt + \sigma^{p^F}(t,P^F_t) dW^F_t,\quad
    P^F_0 = p_F>0, \\
    d P^B_t &= r P^B_t dt + \sigma_B P^B_t dW^B_t \eqqcolon b^{p^B}(t,P^B_t)dt + \sigma^{p^B}(t,P^B_t) dW^B_t,\quad
    P^B_0 = p_B>0, \\
\end{align*}
where $r>0$ is the risk-free interest rate, $\sigma_F>0$ and $\sigma_B>0$ are the volatilities of the feed and biomass prices, respectively. The initial prices of feed and biomass are given by $p_F>0$ and $p_B>0$. The processes $W^F_t$ and $W^B_t$ are two standard independent Brownian motions.

\paragraph*{Farm value.}\label{sec:FarmValue}

To incorporate feeding costs we will consider the following instantaneous cost function:
\begin{align*}
    k(t,w,h,p^F,p^B,u) \coloneqq - h \cdot u \cdot p^F,
\end{align*}
i.e., the cost at time $t\geq 0$ is given by the number of fish $h$, multiplied by the feeding rate $u$ and the feed price $p^F$.
Similarly, we will consider the following terminal reward function at the time of harvesting:
\begin{align*}
    g(t,w,h,p^F,p^B) \coloneqq w \cdot h \cdot p^B,
\end{align*}
where the revenue is given by the biomass $w\cdot h$ multiplied by the biomass price $p^B$.

The value of the farm is the given by the following functional:
\begin{align*}
    J(t,x; u,\tau) \coloneqq 
    \E^{\Q}_{t,x}\left[ 
        \int_t^\tau e^{-r(s-t)} k(s,w_s^{f,u},h_s^{f,u},P^F_s,P^B_s,u_s) ds + 
        e^{-r(\tau-t)} g(\tau,w_\tau^{f,u},h_\tau^{f,u},P^F_\tau,P^B_\tau) 
    \right],
\end{align*}
where $x\coloneqq (w,h,p^F,p^B)$ and $\tau$ is a stopping time representing the time of harvesting and $u$ is the feeding strategy.

We are therefore interested in the following joint optimal control and stopping problem:

\begin{align}
    \sup_{u\in\mathcal{A}, \tau\in\mathcal{T}_{[0,T]}} J(0,x; u,\tau),
    \tag{JCtrlOS}
    \label{eq:JCtrlOS}
\end{align}
where $\mathcal{A}$ is the set of admissible feeding strategies and $\mathcal{T}_{[0,T]}$ the set of stopping times with values in $[0,T]$ for some fixed time horizon $T>0$.

Since both the cost and reward function are dependent on $h$ and $h_t$ is an exponential function, it is equivalent to consider $h_0=1$ and multiply the farm value by the original $h_0$ afterwards. This makes discretizing the state space easier for the numerical application.

Since the feeding strategy $f$ is assumed to be fixed, we will drop in the following the dependence $f$ whenever there is no confusion.

\subsection{HJB Variational Inequality and Free Boundary Problem}\label{sec:HJBVI}
Following \citet[p.~19, Equation 4.2]{Dumitrescu2016}, we can characterize the value function of \eqref{eq:JCtrlOS} under some regularity conditions as the unique viscosity solution of the following local HJB variational inequality  ($x=[w,h,p^F,p^B]$):
\begin{align}
    \max\left\{ \frac{\partial V}{\partial t} + \sup_{u\in[0,\bar{u}]} \mathcal{L}^u V - r V + k(t,x,u), g(t,x) - V \right\} = 0,
    \tag{HJB-VI}
    \label{eq:HJB-VI}
\end{align}
where the infinitesimal generator $\mathcal{L}^u$ is given by
\begin{align*}
    \mathcal{L}^u V(t,x) &\coloneqq 
    \sum_{i\in\{w,h,p^F,p^B\}} b^i(t,x;u) \frac{\partial V}{\partial i}(t,x) + \frac{1}{2} \sum_{i\in\{p^F,p^B\}} (\sigma^i(t,x))^2 \frac{\partial^2 V}{\partial i^2}(t,x) \\&=
    \frac{\partial V}{\partial w}(t,x) \left(\gamma - \gamma^F (f_t-u)^2\right) \left( 1 - \left(\frac{w}{w_\infty}\right)^\nu \right) \\
    &\quad + \frac{\partial V}{\partial h}(t,x) \left(-\mu - \mu^F (f_t-u)^2\right) h \\
    &\quad + \frac{\partial V}{\partial p^F}(t,x) r p^F + \frac{1}{2} \frac{\partial^2 V}{\partial (p^F)^2}(t,x) (\sigma_F p^F)^2 \\
    &\quad + \frac{\partial V}{\partial p^B}(t,x) r p^B + \frac{1}{2} \frac{\partial^2 V}{\partial (p^B)^2}(t,x) (\sigma_B p^B)^2.
\end{align*}
We notice that the control only acts on the weight and population dynamics, as well as the cost function, but not on the price dynamics. Therefore, we can find an analytical expression in terms of the value function $V$ for the optimal control by considering the following first order condition:
\begin{align*}
    \partial_u \left( \mathcal{L}^u V(t,x) + k(t,x,u) \right) = 0.
\end{align*}
This is equivalent to 
\begin{align}
    u(t,[w,h,p^F,p^B]) = f_t-\frac{h p^B}{2 h V_h(t,x) \mu _F-2 V_w(t,x) w \gamma _F \left(\left(\frac{w}{w_{\infty }}\right)^{\nu }-1\right)}.
    \label{eq:FeedbackControl}
\end{align}
Now, this critical point is a local maximum if the second order condition
\begin{align*}
    -2 h V_h \mu _F-2 V_w w \gamma _F \left(1-\left(\frac{w}{w_{\infty }}\right)^{\nu }\right) < 0
\end{align*}
is satisfied.
In our case $h,w,\gamma_F,\mu_F >0$ and by definition of the logistic growth model we have $1-\left(\frac{w}{w_{\infty }}\right)^{\nu }>0$. Therefore, the second order condition is satisfied if $V_h,V_w>0$, which is economically reasonable since the value of the farm should increase with increasing weight and population size.

Since we will restrict the control to be in the interval $[0,\bar{u}]$ for some $\bar{u}>0$, we have to argue for the boundary cases as well. 
For $u=0$, we would have maximal penalty for underfeeding and therefore the optimal control would never be $0$ unless the feed price is $\infty$ as well. If we choose $\bar{u}$ sufficiently large, we can assume that the optimal control will never be greater than $f_t$, since overfeeding is penalized by the additional feeding costs and a decrease in growth and population size and therefore never be optimal in our setting. Hence the optimal control will be given by \eqref{eq:FeedbackControl} projected onto the interval $[0,f_t]$.

Another way to interprete \eqref{eq:HJB-VI}, is in terms of a free boundary problem. We define the \emph{continuation region} by
\begin{align*}
    \mathcal{C} \coloneqq \{ (t,x)\in[0,T]\times\R^4_+ : V(t,x) > g(t,x) \},
\end{align*}
i.e., the region where it is optimal to continue the operation of the farm, and the \emph{stopping region} by
\begin{align*}
    \mathcal{S} \coloneqq \{ (t,x)\in[0,T]\times\R^4_+ : V(t,x) \leq g(t,x) \},
\end{align*}
This allows us to adapt the results by \citet{Sirignano2018} and use PINNs to solve the HJB equation in the continuation region 
\begin{align}
    \left\{
    \begin{aligned}
        \frac{\partial V}{\partial t} + \sup_{u\in[0,\bar{u}]} \mathcal{L}^u V - r V + k(t,x,u) &= 0,&& &(t,x)&\in\mathcal{C}\\
        V(t,x) &= g(t,x),&&  &(t,x)&\in\partial\mathcal{C}\\
        V(T,x) &= g(T,x),&&  &x&\in\R^4_+ 
    \end{aligned}
    \right.
    \tag{HJB-FB}
    \label{eq:HJB-FB}
\end{align}
For our adaption, we will first use \eqref{eq:FeedbackControl} to express the optimal control in terms of the value function and then solve \eqref{eq:HJB-FB}. After that, we will show that a second neural network can be used to approximate the optimal control directly without knowing the feedback form. 
For the stopping region to be learned correctly, we will use an idea similar to \citet{Reppen2025} to make the stopping region a bit fuzzy during the training of the neural network, which improves the accuracy of the stopping region significantly.

\section{Numerical Experiments}\label{sec:NumExp}

For each of the following numerical experiments, we will determine the optimal control and stopping times. After that, we will simulate $8192$ Monte Carlo paths of the underlying Brownian motions and evaluate the farm value for each of the paths using the same realizations of the Brownian motions to make it comparable. This means that the farm value itself will be subject to a Monte-Carlo error, but the comparison of the methods is pathwise and therefore not affected by Monte-Carlo errors. Moreover, we will assume that the feeding rate is given by $f_t\coloneqq f_0 + \eta t$ in all experiments of this section. We chose $\eta$, such that the feeding rate will be in the interval $[0,1]$ to suggest a feeding strategy in percentage. We refer to \Cref{sec:Appendix} for experiments with different feeding strategies, which will show similar results.

In \Cref{tab:Parameters} we summarize the parameters used in the numerical experiments.
\begin{table}[H]
    \centering
    \caption{Parameters used in the numerical experiments.}
    \label{tab:Parameters}
    \begin{tabular}{@{}l l l@{}}
        Parameter & Description & Value \\
        \midrule
        $h_0$          &  Initial population &   $1.000$ \\
        $w_0$          &  Initial weight &   $0.010$ \\
        $p^F_0$        &  Initial feed price &   $0.075$ \\
        $p^B_0$        &  Initial biomass price &   $0.100$ \\
        $\mu_0$        &  Intrinsic mortality rate &   $0.100$ \\
        $\mu_F$        &  Mortality penalty due to under- or overfeeding &   $3.000$ \\
        $\gamma_0$     &  Intrinsic growth rate &   $5.000$ \\
        $\gamma_F$     &  Growth rate penalty due to under- or overfeeding &   $10.000$ \\
        $w_{\infty}$  &  Asymptotic weight &   $3.000$ \\
        $\nu$          & Shape parameter & $0.750$ \\
        $r$            &   Risk-free interest rate & $0.010$ \\
        $\sigma^{p^F}$ &  Volatility of feed price & $0.250$ \\
        $\sigma^{p^B}$ &  Volatility of biomass price & $0.100$ \\
        $T$            &  Time horizon & $3.000$ \\
        $f_0$          &  Initial feed amount & $0.100$ \\
        $\eta$         &  Linear feeding rate & $0.300$ \\
    \end{tabular}
\end{table}

This section is organized as follows. In \Cref{sec:Benchmarks} we will introduce two benchmarks as a sanity check. In \Cref{sec:FDS} we will present a finite difference scheme to solve \eqref{eq:HJB-VI} and in \Cref{sec:DS} we will present a PINN approach to solve \eqref{eq:HJB-FB}. 
\subsection{Benchmarks}\label{sec:Benchmarks}
To judge the improvement to common practice in aquaculture operations and quality of our numerical schemes, we will compare the results to two benchmarks. The first benchmark is the biological optimal feeding strategy $u_t \equiv f_t$ and the second benchmark is a constant feeding strategy $u\equiv 0$. 
Moreover, we will consider for each of the benchmarks the case with and without the option to stop the operation of the farm at any time. For this we set 
\begin{align*}
    \tau^0 \coloneqq T, \quad  \tau^1 \coloneqq \argmax_{ t \in [0,T]} \{ h^{f_t,f_t}_t \cdot w^{f_t,f_t}_t \}, \quad \tau^2 \coloneqq 
    \argmax_{ \tau \in \mathcal{T} } J(0,x;u,\tau).
\end{align*}
The latter, we will determine by the Deep Optimal Stopping (DeepOS) algorithm by \citet{Becker2020}. The results are collected in \Cref{tab:Benchmarks}.
It consists of a matrix for each of the different stopping times and feeding strategies. The values in the table represent the mean farm values over $8192$ Monte Carlo simulations.

\begin{table}[H]
    \caption{Benchmarks used in the numerical experiments. Values in table representing mean farm values over $8192$ Monte Carlo simulations.}
    \label{tab:Benchmarks}
    \centering
    \begin{tabularx}{\linewidth}{@{}l|C|C@{}}
    \diagbox{Stopping Time}{Feeding Strategy}   & $u_t\equiv 0$ & $u_t\equiv f_t$ \\
    \hline
    $\tau^0$ ($=T=3.0$)                             & 0.001  & 0.1175  \\
    $\tau^1$ ($\approx 2.176$)                           & 0.0285  & 0.1732  \\
    $\tau^2$                      & 0.0565 \newline($\E[\tau^2] \approx 1.533$) & 0.1798\newline($\E[\tau^2] \approx 1.947$)  \\
    \end{tabularx}
\end{table}

We can conclude that our numerical schemes for the joint optimal control and stopping problem should at least outperform the best benchmark, which is the biological optimal feeding strategy with optimal stopping.
\subsection{Finite Difference Solver}\label{sec:FDS}
In order to solve \eqref{eq:HJB-VI} numerically, we will use the dynamic programming principle and iterate backwards in time. 
For this, we will combine the methods presented by \citet[pp.~142\,ff., 327\,ff.]{Kushner2001} for the optimal stopping and control problem.

To this end, we will discretize the state space and time interval and use finite differences to approximate the spatial derivatives. We will assume homogeneous grids for all dimensions and denote their step sizes by $\Delta t, \Delta w, \Delta h, \Delta p^F, \Delta p^B >0$. 
For the time derivative, we will use an explicit Euler scheme, for the first order spatial derivatives, we will use an upwind scheme and for the second order spatial derivatives, we will use a central difference scheme. To be more precise, we will approximate the derivatives as follows:
\begin{align*}
    \partial_t f(t,x) & \approx \frac{f(t+\Delta t,x)-f(t,x)}{\Delta t}, \\
    \partial_x f(t,x) & \approx 
    \begin{cases}
        \frac{f(t,x+\Delta x)-f(t,x)}{\Delta x}, & \text{if } b(x,u)\geq 0, \\
        \frac{f(t,x)-f(t,x-\Delta x)}{\Delta x}, & \text{if } b(x,u)< 0,
    \end{cases}\\
    \partial_{xx} f(t,x) &\approx \frac{f(t,x+\Delta x)-2 f(t,x)+f(t,x-\Delta x)}{(\Delta x)^2}.
\end{align*}
For the upwind scheme, we will use the following shorthand notation:
\begin{align*}
    f^+(t,x) = \max(f(t,x),0), \quad f^-(t,x) = \max(-f(t,x),0), \quad \abs{f(t,x)} = f^+(t,x) + f^-(t,x).
\end{align*}

This leads to the following explicit scheme (suppressing the arguments of $b^i$ and $\sigma^i$ for better readability):
\begin{align*}
    \hspace*{1em}&\hspace*{-1em}
    \tilde{V}^u(t_n,w_i,h_j,p^F_k,p^B_l) \\&\coloneqq 
    V(t_n,w_i,h_j,p^F_k,p^B_l) \left(1- \Delta t \sum_{i\in \{w,h,p^F,p^B\}} \frac{\abs{b^i}}{\Delta i} - \Delta t \sum_{i\in \{p^F,p^B\}} \frac{(\sigma^i)^2}{(\Delta i)^2}\right) \\& \quad+
    V(t_n,w_i+\Delta w,h_j,p^F_k,p^B_l) \frac{\Delta t}{\Delta w} (b^w)^+ \\& \quad+
    V(t_n,w_i-\Delta w,h_j,p^F_k,p^B_l) \frac{\Delta t}{\Delta w} (b^w)^- \\& \quad+
    V(t_n,w_i,h_j+\Delta h,p^F_k,p^B_l) \frac{\Delta t}{\Delta h} (b^h)^+ \\& \quad+
    V(t_n,w_i,h_j-\Delta h,p^F_k,p^B_l) \frac{\Delta t}{\Delta h} (b^h)^- \\& \quad+
    V(t_n,w_i,h_j,p^F_k+\Delta p^F,p^B_l) \left(\frac{\Delta t}{\Delta p^F} (b^{p^F})^+ + \frac{\Delta t}{2 (\Delta p^F)^2} (\sigma^{p^F})^2\right) \\& \quad+
    V(t_n,w_i,h_j,p^F_k-\Delta p^F,p^B_l) \left(\frac{\Delta t}{\Delta p^F} (b^{p^F})^- + \frac{\Delta t}{2 (\Delta p^F)^2} (\sigma^{p^F})^2\right) \\& \quad+
    V(t_n,w_i,h_j,p^F_k,p^B_l+\Delta p^B) \left(\frac{\Delta t}{\Delta p^B} (b^{p^B})^+ + \frac{\Delta t}{2 (\Delta p^B)^2} (\sigma^{p^B})^2\right) \\& \quad+
    V(t_n,w_i,h_j,p^F_k,p^B_l-\Delta p^B) \left(\frac{\Delta t}{\Delta p^B} (b^{p^B})^- + \frac{\Delta t}{2 (\Delta p^B)^2} (\sigma^{p^B})^2\right) 
\end{align*}

\begin{align*}
    &V(T,w_i,h_j,p^F_k,p^B_l) = g(T,w_i,h_j,p^F_k,p^B_l), \\
    &V(t_n,w_i,h_j,p^F_k,p^B_l) =\\&\quad \max\Bigg\{ g(t_n,w_i,h_j,p^F_k,p^B_l),\sup_u \left(1 - r \Delta t\right) \tilde{V}(t_n,w_i,h_j,p^F_k,p^B_l)+k(t_{n},w_i,h_j,p^F_k,p^B_l)\Delta t\Bigg\} 
\end{align*}
If one removes the maximum between the continuation value and exercise value $g$, one obtains the scheme for only the optimal control problem without stopping.

In our numerical implementation, we will discretize the control space as well and take the maximum over the discrete set of controls. Then we use linear interpolation for both the value function as well as the control function, which is determined in each time step by saving the maximum over the the grid $u$ for each point $x=[w,h,p^F,p^B]$. This will serve as a benchmark for the deep learning approach presented in the next section. This allows us to parallelize the computation of $\tilde{V}^u$ for all controls $u$ and spatial grid points and then take the maximum afterwards, which makes it ideal for GPU computing. The dimensionality of this problem is the limit for the finite difference approach using this method, since it becomes too memory intensive for higher dimensions. 

For the numerical implementation, we used homogeneous grids for all dimensions and a summary of the grid sizes is given in \Cref{tab:FDGridSizes}.
\begin{table}[H]
    \caption{Grid sizes used in the finite difference solver.}
    \label{tab:FDGridSizes}
    \centering
    \begin{tabular}{@{}l l l l@{}}
        Dimension & Lower Bound & Upper Bound & Grid Size \\
        \midrule
        Time $t$ & 0 & T & 2048 \\
        Weight $w$ & $\frac{w_0}{2}$ & $1.1\, w_\infty$ & 64 \\
        Population $h$ & $\frac{h_0}{10}$ & $1.1\, h_0 $ & 64 \\
        Feed Price $p^F$ & 0.0019 & 0.3856 & 32 \\
        Biomass Price $p^B$ & 0.0055 & 0.2635 & 32 \\
        Control $u$ & 0 & 1 & 64 \\
    \end{tabular}
\end{table}
The grids were chosen such that the memory consumption fits on a single GPU with $24$GB of memory to speed up the computation. The time step was chosen such that the explicit scheme is stable. The spatial grid boundaries were chosen such that they do not influence the solution too much and the feed price and biomass price boundaries were determined by first simulating the processes and taking the minimum and maximum of the samples. The control grid was chosen such that the optimal control can be approximated well enough.

Let us compare the following scenarios:
\begin{compactenum}
    \item Only Optimal Control (denoted by $\hat{u}$) and no stopping, i.e. $\hat{u} = \argmax_u J(0,x,u,T)$;
    \item Only Optimal Control and stopping at $\tau^1$, i.e. $J(0,x,\hat{u},\tau^1)$;
    \item Only Optimal Control and then optimal stopping at $\tau^2$ (using DeepOS), i.e. $\tau^2 = \argmax_\tau J(0,x,\hat{u},\tau)$;
    \item Only Optimal Control till time $\tau^1$ (denoted by $\tilde{u}$) and no stopping, i.e. $J(0,x,\tilde{u},\tau^1)$;
    \item Only Optimal Control and then approximate stopping at $\tau^3$ (using DeepOS), i.e. $\tau^3 = \argmax_\tau J(0,x,\tilde{u},\tau)$;
    \item Joint Optimal Control and Stopping, i.e. $u^\ast, \tau^\ast = \argmax_{u,\tau} J(0,x,u,\tau)$.
\end{compactenum}
The results of these scenarios are summarized in \Cref{tab:FDResults}, in which each row corresponds to one of the scenarios and the second row gives the expected stopping time $\E[\tau]$ and the third row the mean farm value $J(0,x;u,\tau)$ over $8192$ Monte Carlo simulations.

\begin{table}[H]
    \caption{Results of the finite difference solver. Values in table representing mean farm values over $8192$ Monte Carlo simulations.}
    \label{tab:FDResults}
    \centering
    \begin{tabular}{@{}lcc@{}}
        Scenario & $\E[\tau]$ & $J(0,x;u,\tau)$ \\
        \midrule
        1. $\hat{u}, \tau^0$ & 3.0 & 0.1259\\
        2. $\hat{u}, \tau^1$ & 2.1764 & 0.1760\\
        3. $\hat{u}, \tau^2$ & 1.9611 & 0.1793\\
        4. $\tilde{u}, \tau^1$ & 2.1764 & 0.1775\\
        5. $\tilde{u}, \tau^3$ & 1.9410 & 0.1826\\
        6. $u^\ast, \tau^\ast$ & 1.9652 & 0.1841\\
    \end{tabular}
\end{table}

We can see that the joint optimal control and stopping problem leads to the highest farm value, as expected. Moreover, we can see that the optimal control till time $\tau^1$ leads to a higher farm value than the optimal control till time $T$, which shows that it is beneficial to stop the operation of the farm at some point in time. It is noteworthy that the optimal stopping time from the joint optimal control and stopping problem is later than the optimal stopping time when using the optimal control till time $\tau^1$. 

\begin{figure}
    \centering
    \includegraphics[width=.9\columnwidth]{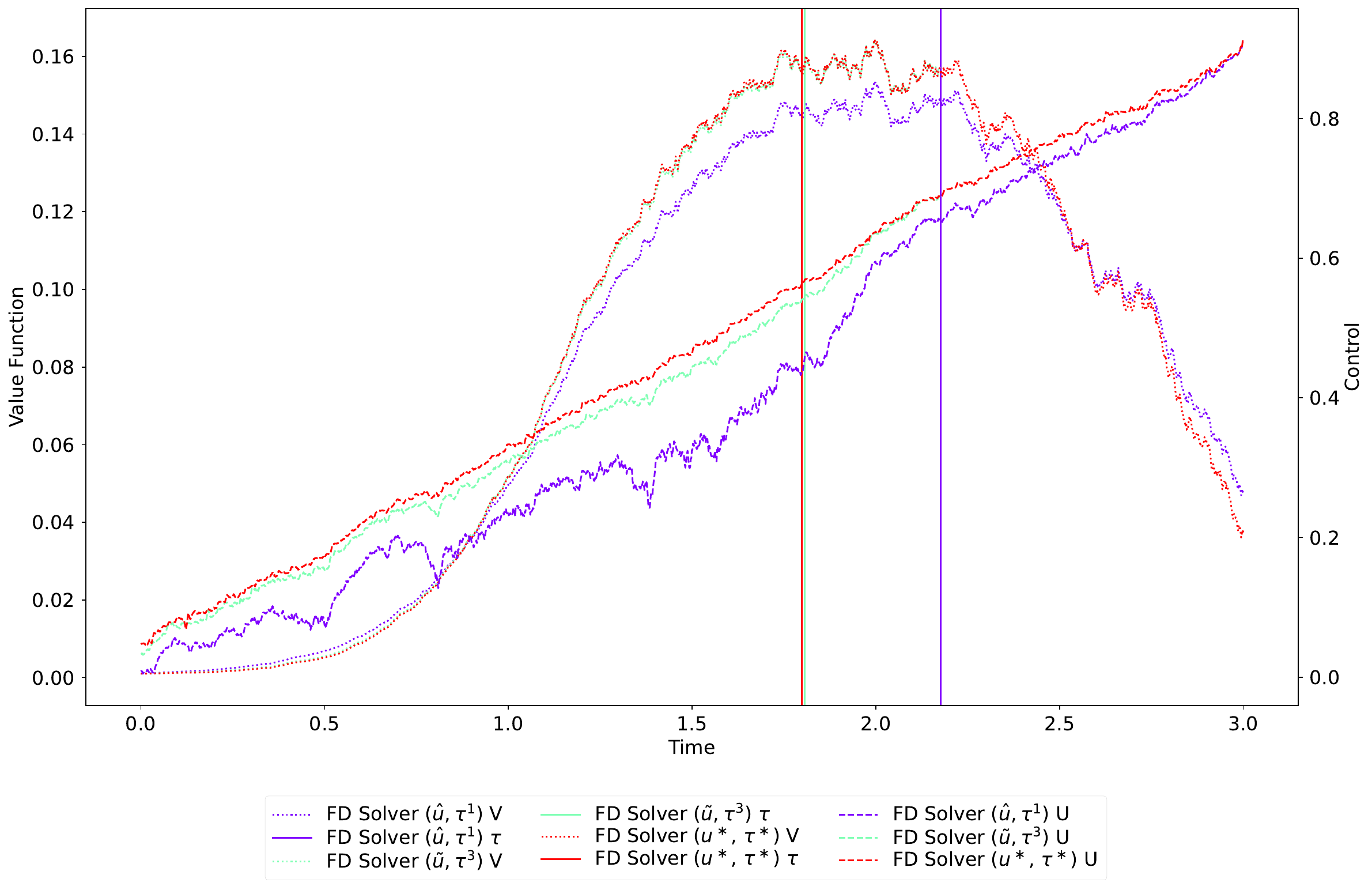}
    \caption{Comparison of the optimal control and stopping strategy obtained by the finite difference solver. The plots represent one trajectory of the simulated farm values.}
    \label{fig:FD_comparison}
\end{figure}

In \Cref{fig:FD_comparison} we compare the optimal control and stopping strategy obtained by the finite difference solver. The plots represent one trajectory of the simulated farm values. Bold lines represent the stopping time, dashed lines the control and dotted lines the farm value over time. 
The left y-axis represents the farm value and the right y-axis the control value. The color purple corresponds to scenario 2, green to scenario 5 and red to scenario 6 in \Cref{tab:FDResults}.

We can see that $\hat{u}$ is the most different from the other two controls $\tilde{u}$ and $u^\ast$, due to the fact that the farm needs to be controlled till $T$, where the farm value of scenario 2 is higher than the one of scenario 6, which is expected. For scenario 5 and 6, the controls and farm value are quite similar but the later harvest in scenario 6 leads to a slightly higher farm value. 


\subsection{Deep Solver}\label{sec:DS}
In this section, we will present a deep learning approach to solve \eqref{eq:HJB-FB} using Physics-Informed Neural Networks (PINNs) as introduced by \citet{Raissi2019} and adapted to free boundary problems by \citet{Sirignano2018}. 

\paragraph*{Neural Network Architecture for the Value Function.}
We will use a fully connected feedforward neural network with $3$ hidden layers and $32$ neurons per layer and use the same bounds as for the finite difference solver (see \Cref{tab:FDGridSizes}) for sampling grid points for the PDE. As in \citet{Sirignano2018}, the activation function will be the hyperbolic tangent function. The input of the neural network will be the time $t$ and the state variables $x=[w,h,p^F,p^B]$ and the output will be the value function $V(t,x)$. Since this is a low-dimensional setting from a Deep Learning perspective, we will can compute the derivatives of the neural network using automatic differentiation without performance issues.
The loss function (see \citet{Sirignano2018}) will be a combination of the PDE loss, the boundary loss and the terminal condition loss. To be more precise, we will use the following loss functions $L_{\text{PDE}},L_{\text{FB}} ,L_{\text{T}}$:
\begin{align*}
    L &\coloneqq \mathrm{MSE}(L^u_{\text{PDE}}) + \mathrm{MSE}(L_{\text{FB}}) + \mathrm{MSE}(L_{\text{T}}), \\
    L^u_{\text{PDE}} &\coloneqq \frac{\partial V}{\partial t} + \mathcal{L}^u V - r V + k(t,x,u), \quad \text{ on } V(t,x) > g(t,x), \\ 
    L_{\text{FB}} &\coloneqq V(t,x) - g(t,x), \quad \text{ on } V(t,x) \leq g(t,x), \\
    L_{\text{T}} &\coloneqq V(T,x) - g(T,x).
\end{align*}
The acronym $\mathrm{MSE}$ stands for the mean squared error. The loss function $L_{\text{PDE}}$ enforces the PDE to be satisfied in the continuation region, $L_{\text{FB}}$ enforces the value function to be greater or equal to the exercise value on the free boundary and $L_{\text{T}}$ enforces the terminal condition to be satisfied.

\paragraph*{Approximation of Control Function.}
To approximate the optimal control, we will compare three different approaches. We will denote the output of the three approaches after optimization by $u^1,u^2,u^3$, respectively.

The first approach is to use \eqref{eq:FeedbackControl} to express the optimal control in terms of the value function and its derivatives. The other two approaches aim at situations, where the feedback form is not known, to develop a black-box approach to approximate the optimal control directly. 

The second approach is to use a second neural network with the same architecture as the value function neural network to approximate the optimal control directly. We use ReLU activation functions for this network and use the absolute value on the output of the network to ensure positive controls. The input of this neural network will be again $t$ and $x$ and the output will be the control $u(t,x)$. The loss function for this approach will be the following:
\begin{align*}
  L_{\text{ctrl}} = \mathrm{ME}(-\mathcal{L}^u V + k(t,x,u)),
\end{align*}
where $\mathrm{ME}$ stands for the mean error. The idea behind this loss function is that we want to maximize the infinitesimal generator $\mathcal{L}^u V$ plus the cost term $k(t,x,u)$ with respect to $u$ and therefore minimize its negative.

For the third approach, we use the exact same neural network architecture as for approach two but change the loss function. The loss function for this approach will be the following: First we evaluate the infinitesimal $\mathcal{L}^u V + k(t,x,u)$ on a predefined grid of controls $u_1,\ldots,u_m$ and determine its maximum, denoted by $\hat{y}$. Then, we compute controls $\hat{u}$ using approach two and evaluate $\mathcal{L}^{\hat{u}} V + k(t,x,\hat{u})$, denoted by $y$. The loss function will then be the following:
\begin{align*}
 L_{\text{ctrl}} = - (\hat{y} - y)^+.
\end{align*}
This ensures the control network outperforms the optimal control on the discrete grid of controls.

\paragraph*{Data Generation.}
To train the neural networks, we need to generate a dataset of grid points in the state space and time. We will use the same bounds as for the finite difference solver (see \Cref{tab:FDGridSizes}) and sample the points uniformly in the time interval and state space for each training epoch. In each training step we will \emph{balance} the number of points in the continuation region and on the free boundary by using rejection sampling. For this, we first sample uniformly in the whole domain and then keep $4096$ points where $V(t,x) > g(t,x)$ according to the current neural network approximation. For sampling points on the free boundary, we will use a similar approach as \citet{Reppen2025} and sample uniformly in the whole domain and then keep $4096$ points where $\abs{V(t,x) - g(t,x)} < \varepsilon \cdot u$, where $u\sim \mathrm{Uniform}(0,1)$ are uniform random numbers in the interval $0$ and $1$ for some small $\varepsilon>0$. The uniform random numbers ensure that the network does not learn the shift by $\varepsilon$. This makes the stopping region a bit fuzzy during training, which improves the accuracy of the stopping region significantly.

\paragraph*{Training.}
The training of the neural networks will be performed using the Adam optimizer with a learning rate given by the following schedule 
\begin{align*}
    \mathrm{lr}(e) = \mathrm{lr}_0 \cdot \left(\frac{1}{2}\right)^{\lfloor \frac{e}{1000} \rfloor},
\end{align*}
where $\mathrm{lr}_0 = 5\cdot 10^{-3}$ is the initial learning rate and $e$ is the current epoch.
We will use batches of size $4096$ and train the value function networks for $10000$ epochs. For each epoch of the value function network, we train the control networks $5$ times with a fixed learning rate $5 \cdot 10^{-4}$ also with Adam optimizer. We chose $\varepsilon=0.01$ for the fuzzy boundary sampling, which seems to perform best in our experiments.

\paragraph*{Results.}
We will compare the three approaches for approximating the optimal control presented above to the solution obtained by the finite difference approach. The results are summarized in \Cref{tab:DLResults}, in which each row corresponds to one of the Deep Learning approaches and the second row gives the expected stopping time $\E[\tau^\varepsilon]$ and the third row the mean farm value $J(0,x;u,\tau)$ over $8192$ Monte Carlo simulations. The last column shows the difference to the finite difference solution. The optimal stopping time is determined in the simulation by the first time the value function is less or equal to the exercise value plus a small tolerance $\varepsilon=0.01$.
\begin{table}[H]
    \caption{Results of the deep learning solver. Values in table representing mean farm values over $8192$ Monte Carlo simulations.}
    \label{tab:DLResults}
    \centering
    \begin{tabular}{@{}lccc@{}}
        Method & $\E[\tau^\varepsilon]$ & $J(0,x;u^i,\tau^\varepsilon)$ & $\abs{J(0,x;u^i,\tau^\varepsilon) - J(0,x;u^\ast,\tau^\ast)}$ \\
        \midrule
        1. Feedback Control $u^1$ & 1.8729 & 0.1835 & 0.0005 \\
        2. Control Network $u^2$ & 1.9718 & 0.1835 & 0.0005 \\
        3. Control Network $u^3$ & 2.1604 & 0.1784 & 0.0054 \\
    \end{tabular}
\end{table}
We see that all three approaches lead to comparable results, with approach 2 leading to the highest farm value. However, all three approaches lead to slightly lower farm values than the finite difference approach, which suggests that the either the optimal stopping time is not determined correctly by the PINNs or the control network needs further refinement.

\begin{figure}
    \centering
    \includegraphics[width=.9\columnwidth]{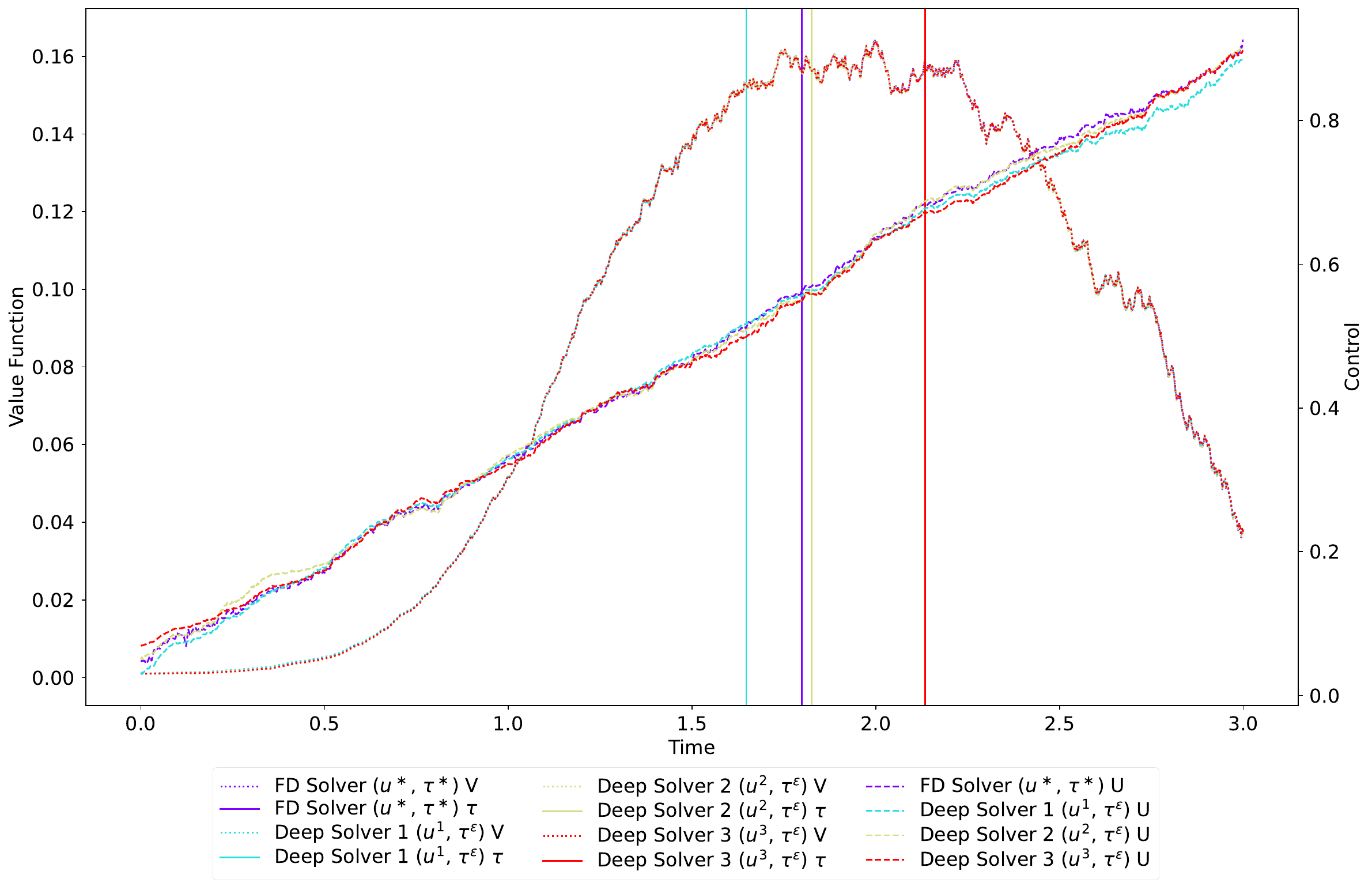}
    \caption{Comparison of the optimal control and stopping strategy obtained by the deep learning solvers and finite difference solver. The plots represent one trajectory of the simulated farm values.}
    \label{fig:DL_comparison}
\end{figure}

In \Cref{fig:DL_comparison}, we compare the optimal control and stopping strategy obtained by the deep learning solver using all three approaches. The plots represent one trajectory of the simulated farm values. Bold lines represent the stopping time, dashed lines the control and dotted lines the farm value over time. The left y-axis represents the farm value and the right y-axis the control value. The color purple corresponds to the finite difference solution, blue to approach 1, yellow to approach 2 and red to approach 3 in \Cref{tab:DLResults}.

We can see that the controls all very similar. To investigate, if the optimal stopping time is determined correctly, we compute the optimal stopping time with the DeepOS algorithm by \citet{Becker2020} using the control $u^1,u^2,u^3$. The results are summarized in \Cref{tab:DLResultsOS}, in which each row corresponds to one of the Deep Learning approaches and the second row gives the expected stopping time $\E[\tau^{\text{OS}}]$ and the third row the mean farm value $J(0,x;u,\tau^{\text{OS}})$ over $8192$ Monte Carlo simulations. The last column shows the difference to the finite difference solution.
\begin{table}[H]
    \caption{Results of the deep learning solver combined with DeepOS. Values in table representing mean farm values over $8192$ Monte Carlo simulations.}
    \label{tab:DLResultsOS}
    \centering
    \begin{tabular}{@{}lccc@{}}
        Method & $\E[\tau^\text{OS}]$ & $J(0,x;u^i,\tau^\text{OS})$ & $\abs{J(0,x;u^i,\tau^\text{OS}) - J(0,x;u^\ast,\tau^\ast)}$ \\
        \midrule
        1. Feedback Control $u^1$ & $1.9597$ & $0.1833$ & $0.0007$ \\
        2. Control Network $u^2$ & $1.9654$ & $0.1825$ & $0.0015$ \\
        3. Control Network $u^3$ & $1.9640$ & $0.1829$ & $0.0011$ \\
    \end{tabular}
\end{table}
We can see that the combination of PINNs and DeepOS can remedy the short-comings of PINNs for approximating the free boundary, since the mean of the optimal stopping times look now very similar to $\tau^\ast$. The values for all approaches are now lower than without DeepOS, which can be due to a too coarse time grid for the DeepOS model or random chance of the fuzzy stopping rule from before. Thus, the network for the value function needs to be trained more accurately to obtain better results.
The approach with feedback control is still the most accurate of the three. This suggests that the approximation of the optimal control of the other two methods is not as accurate as the feedback form and needs further investigation in the future, either by improving the loss function, the neural network architecture or the training procedure. In the experiments in the appendix \Crefrange{sec:ExponentialFeedingRate}{sec:SinusoidalFeedingRate}, we can see that the control network approaches have a similar performance as in \Cref{tab:DLResultsOS} as long as there are not so many oscillations in the feeding rate. Which is expected, since the control network needs to learn a more complex mapping from the state space to the control space.

\section{Conclusion}\label{sec:Conclusion}
In this work, we presented a joint optimal control and stopping problem for aquaculture operations. We formulated the problem as a variational inequality and a free boundary problem and presented two numerical approaches to solve it. The first approach is a finite difference scheme, which can be used as a benchmark for the second approach, which is a deep learning approach using PINNs (combined with DeepOS). Both methods took less than 10 minutes to compute the solution on a single GPU.

In our tests, we could see that the finite difference approach is still working in this 5-dimensional setting. We showed that the highest farm values are obtained by the joint optimal control and stopping problem and neglecting either the optimal control or stopping leads to lower farm values.

We demonstrated that PINNs can be used to solve the joint optimal control and stopping problem as well. We presented three different approaches to approximate the optimal control, one using the feedback form and two using a second neural network to approximate the optimal control directly. We found that the approximation of the stopping time is not as accurate as the finite difference approach, which results in lower farm values. However, we could see that a combination of PINNs and DeepOS can remedy the short-comings of PINNs for approximating the free boundary.

The results of the deep learning approach are comparable to the finite difference approach, using either the feedback controls or a second neural network to approximate the optimal control directly. This shows that the presented deep learning approach is a viable alternative to the finite difference approach, and opens the door to numerical methods in higher dimensions, where the finite difference approach becomes infeasible due to the curse of dimensionality. That being said, a careful choice for the control network architecture, loss function and training procedure is necessary to obtain good results.

The methodology can be easily adapted to other joint optimal control and stopping problems in finance and economics. Future research could focus on improving the approximation of the stopping time by PINNs and investigating the performance of the deep learning approach in higher dimensions.

\printbibliography

@article{EK2025,
title = {On the impact of biological risk in aquaculture valuation and decision making},
journal = {Aquaculture},
volume = {603},
pages = {742368},
year = {2025},
issn = {0044-8486},
doi = {10.1016/j.aquaculture.2025.742368},
author = {Christian Oliver Ewald and Kevin Kamm}
}

@article{EK2024,
  author    = {Christian Oliver Ewald and Kevin Kamm},
  title     = {On the impact of feeding cost risk in aquaculture valuation and decision making},
  journal   = {Quantitative Finance},
  volume    = {24},
  number    = {9},
  pages     = {1341--1352},
  year      = {2024},
  doi       = {10.1080/14697688.2024.2308069}
}

@article{Reisinger2019,
  author       = {Roxana Dumitrescu and Christoph Reisinger and Yufei Zhang},
  title        = {Approximation Schemes for Mixed Optimal Stopping and Control Problems with Nonlinear Expectations and Jumps},
  journal      = {Applied Mathematics \& Optimization},
  volume       = {83},
  number       = {-},  
  pages        = {1387--1429},
  year         = {2019},
  doi          = {10.1007/s00245-019-09591-0}
}

@article{Dumitrescu2016,
  author       = {Roxana-Larisa Dumitrescu and Marie-Claire Quenez and Agnès Sulem},
  title        = {A Weak Dynamic Programming Principle for Combined Optimal Stopping/Stochastic Control with $\mathcal{E}^{f}$-expectations},
  journal      = {SIAM Journal on Control and Optimization},
  volume       = {54},
  number       = {4},
  pages        = {2090--2115},
  year         = {2016},
  doi          = {10.1137/15M1027012}
}

@article{Sirignano2018,
  author    = {Justin Sirignano and Konstantinos Spiliopoulos},
  title     = {DGM: A deep learning algorithm for solving partial differential equations},
  journal   = {Journal of Computational Physics},
  volume    = {375},
  pages     = {1339--1364},
  year      = {2018},
  doi       = {10.1016/j.jcp.2018.08.029}
}

@book{Kushner2001,
  author    = {Harold J. Kushner and Paul Dupuis},
  title     = {Numerical Methods for Stochastic Control Problems in Continuous Time},
  series    = {Stochastic Modelling and Applied Probability},
  volume    = {24},
  edition   = {2},
  publisher = {Springer, New York, NY},
  year      = {2001},
  isbn      = {978-0-387-95139-3, 978-1-4612-6531-3, 978-1-4613-0007-6},
  doi       = {10.1007/978-1-4613-0007-6}
}

@article{Reppen2025,
  author    = {Andres Max Reppen and Halil Mete Soner and Valentin Tissot-Daguette},
  title     = {Neural optimal stopping boundary},
  journal   = {Mathematical Finance},
  volume    = {35},
  number    = {2},
  pages     = {441--469},
  year      = {2025},
  doi       = {10.1111/mafi.12450}
}

@misc{Becker2020,
      title={Deep optimal stopping}, 
      author={Sebastian Becker and Patrick Cheridito and Arnulf Jentzen},
      year={2020},
      eprint={1804.05394},
      archivePrefix={arXiv},
      primaryClass={math.NA},
      url={https://arxiv.org/abs/1804.05394}, 
}

@article{Raissi2019,
title = {Physics-informed neural networks: A deep learning framework for solving forward and inverse problems involving nonlinear partial differential equations},
journal = {Journal of Computational Physics},
volume = {378},
pages = {686-707},
year = {2019},
issn = {0021-9991},
doi = {10.1016/j.jcp.2018.10.045},
author = {M. Raissi and P. Perdikaris and G.E. Karniadakis},
}

\appendix
\section{Appendix}\label{sec:Appendix}
In this appendix, we will present additional numerical experiments for different parameter settings. In particular, we will consider different given biological optimal feeding strategies $f$ to study if the Deep Learning approach can handle more complex feeding strategies as well. Since the code is publicly available, we encourage the reader to try out different parameter settings as well.

\subsection{Exponential Feeding Rate}\label{sec:ExponentialFeedingRate}
In this section, we will investigate the impact of an exponential feeding rate on the optimal control and stopping strategies. The exponential feeding rate is defined as:
\begin{equation*}
f(t) = f_0 e^{\lambda t}
\end{equation*}
where $f_0>0$ is the initial feeding rate and $\lambda>0$ is the growth rate. We keep $f_0$ as in \Cref{tab:Parameters} and choose $\lambda$ such that we normalize the feeding rate over the time horizon $T$, i.e.
\begin{equation*}
    f_0 e^{\lambda T} = 1 \Leftrightarrow \lambda = \frac{\ln(1/f_0)}{T}.
\end{equation*}

\begin{table}[H]
    \caption{Results of the deep learning solver combined with DeepOS compared to the finite difference solver. Values in table representing mean farm values over $8192$ Monte Carlo simulations.}
    \label{tab:ResultsEFR}
    \centering
    \begin{tabular}{@{}lccc@{}}
        Method                                      & $\E[\tau]$    & $J(0,x;u,\tau)$   & $\abs{J(0,x;u^i,\tau^\text{OS}) - J(0,x;u^\ast,\tau^\ast)}$ \\
        \midrule
        Benchmark $u=f_t, \tau^\text{OS}$        & $2.025$      & $0.2003$          & -  \\
        FD solver $u^\ast, \tau^\ast$            & $2.039$      & $0.2045$          & -  \\
        Feedback Control $u^1, \tau^\text{OS}$   & $2.035$      & $0.2036$          & $0.0008$ \\
        Control Network $u^2, \tau^\text{OS}$    & $2.038$      & $0.2034$          & $0.0011$ \\
        Control Network $u^3, \tau^\text{OS}$    & $2.058$      & $0.2031$          & $0.0013$ \\
    \end{tabular}
\end{table}

\subsection{Logistic Feeding Rate}\label{sec:LogisticFeedingRate}
In this section, we will investigate the impact of a logistic feeding rate on the optimal control and stopping strategies. The logistic feeding rate is defined as:

\begin{equation*}
f(t) = f_0 + \frac{L-f_0}{1 + e^{-k(t - t_I)}},
\end{equation*}

where $L>f_0>0$ is the maximum feeding rate, $f_0>0$ is the initial feeding rate, $k>0$ is the growth rate, and $t_I\in (0,T)$ is the inflection point.

We set $L=1$, $k=2.5$ and $t_I = T/2$ for our numerical experiments and keep $f_0$ as in \Cref{tab:Parameters}.

\begin{table}[H]
    \caption{Results of the deep learning solver combined with DeepOS compared to the finite difference solver. Values in table representing mean farm values over $8192$ Monte Carlo simulations.}
    \label{tab:ResultsLFR}
    \centering
    \begin{tabular}{@{}lccc@{}}
        Method                                   & $\E[\tau]$    & $J(0,x;u,\tau)$   & $\abs{J(0,x;u^i,\tau^\text{OS}) - J(0,x;u^\ast,\tau^\ast)}$ \\
        \midrule
        Benchmark $u=f_t, \tau^\text{OS}$        & $1.9268$      & $0.1855$          & -  \\
        FD solver $u^\ast, \tau^\ast$            & $1.9440$      & $0.1896$          & -  \\
        Feedback Control $u^1, \tau^\text{OS}$   & $1.9405$      & $0.1889$          & $0.0007$ \\
        Control Network $u^2, \tau^\text{OS}$    & $1.9471$      & $0.1879$          & $0.0016$ \\
        Control Network $u^3, \tau^\text{OS}$    & $1.9526$      & $0.1882$          & $0.0013$ \\
    \end{tabular}
\end{table}

\subsection{Sinusoidal Feeding Rate}\label{sec:SinusoidalFeedingRate}
In this section, we will investigate the impact of a sinusoidal feeding rate on the optimal control and stopping strategies. The sinusoidal feeding rate is defined as:

\begin{equation*}
f(t) = a \sin\left(\frac{2\pi}{t_p} t \right) + b t + f_0
\end{equation*}

where $a>0$ is the amplitude, $t_p>0$ is the period, $b\geq 0$ is the linear drift and $f_0>0$ is the initial feeding rate. We will choose $a=0.1$, $t_p=T/12$, $\phi=0$ and keep $f_0$ as in \Cref{tab:Parameters}. We will determine $b$ such that the feeding rate is in the interval $[0,1]$ over time, i.e.
we first consider the maximum of the sine function over $[0,T]$, i.e.
\begin{align*}
    \max_{t\in [0,T]} a \sin\left(\frac{2\pi}{t_p} t \right) = a, \text{ for} \quad t = \frac{t_p}{4} + 2\pi n, \quad n \in \mathbb{Z}.
\end{align*}
Since $f(t)$ is in mean increasing, we can substitute the maximum into the normalization condition to obtain
\begin{align*}
    &1 = \tilde{f}(T) \coloneqq a + b T + f_0 \\&
    \Leftrightarrow b = \frac{1 - a - f_0}{T}.
\end{align*}

\begin{table}[H]
    \caption{Results of the deep learning solver combined with DeepOS compared to the finite difference solver. Values in table representing mean farm values over $8192$ Monte Carlo simulations.}
    \label{tab:ResultsSFR}
    \centering
    \begin{tabular}{@{}lccc@{}}
        Method                                   & $\E[\tau]$    & $J(0,x;u,\tau)$   & $\abs{J(0,x;u^i,\tau^\text{OS}) - J(0,x;u^\ast,\tau^\ast)}$ \\
        \midrule
        Benchmark $u=f_t, \tau^\text{OS}$        & $1.9654$       & $0.1840$          & -  \\
        FD solver $u^\ast, \tau^\ast$            & $1.9934$      & $0.1877$          & -  \\
        Feedback Control $u^1, \tau^\text{OS}$   & $1.9745$      & $0.1870$          & $0.0007$ \\
        Control Network $u^2, \tau^\text{OS}$    & $1.9786$      & $0.1809$          & $0.0067$ \\
        Control Network $u^3, \tau^\text{OS}$    & $1.9819$      & $0.1819$          & $0.0057$ \\
    \end{tabular}
\end{table}

\begin{figure}
    \label{fig:DL_comparison_SFR}
    \centering
    \includegraphics[width=.75\textwidth]{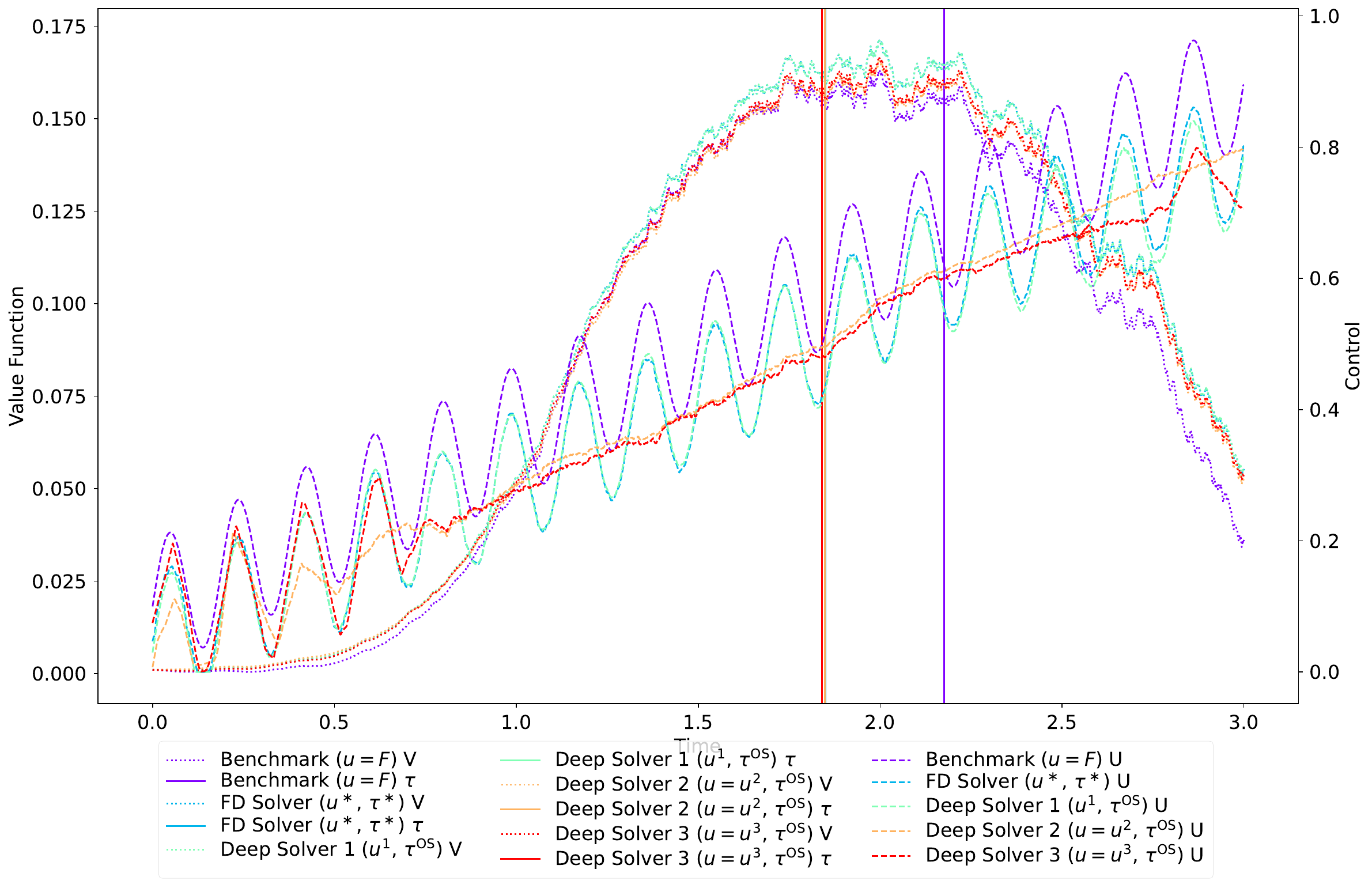}
    \caption{Comparison of the optimal control and stopping strategy obtained by the deep learning solvers and finite difference solver for sinusoidal feeding rate. The plots represent one trajectory of the simulated farm values.}
\end{figure}

In this case, we can see that the control network approaches lead to worse results than the feedback control approach. In \Cref{fig:DL_comparison_SFR}, we compare the optimal control and stopping strategy obtained by the deep learning solver using all three approaches. The plots represent one trajectory of the simulated farm values. Bold lines represent the stopping time, dashed lines the control and dotted lines the farm value over time. The left y-axis represents the farm value and the right y-axis the control value. The color purple corresponds to the benchmark solution with $u=f_t$,
blue to the finite difference solution, green to approach 1, orange to approach 2 and red to approach 3 in \Cref{tab:ResultsSFR}.
We can see that the control networks approximate the sinusoidal feeding rate quite well in the beginning, but deviate more and more from the benchmark solution later on. This suggests that the control networks need further refinement to approximate more complex feeding strategies well. However, the feedback control approach still leads to good results, which shows that the PINN approach is still viable for more complex feeding strategies.
\end{document}